\documentclass[10pt]{notices}

\usepackage{amsmath}
\usepackage{amsfonts}
\usepackage{amssymb}
\usepackage{url}
\usepackage{amsthm}
\usepackage{amssymb,amsmath}

\usepackage[normalem]{ulem}
\usepackage{xcolor}
\hypersetup{
    colorlinks=true,
    linkcolor=blue,
    citecolor=blue,
    urlcolor=blue
}
\def\mvint_#1{\mathchoice
          {\mathop{\vrule width 6pt height 3 pt depth -2.5pt
                  \kern -9pt \intop}\limits_{\kern -3pt #1}}%
          {\mathop{\vrule width 5pt height 3 pt depth -2.6pt
                  \kern -6pt \intop}\nolimits_{#1}}%
          {\mathop{\vrule width 5pt height 3 pt depth -2.6pt
                  \kern -6pt \intop}\nolimits_{#1}}%
          {\mathop{\vrule width 5pt height 3 pt depth -2.6pt
                  \kern -6pt \intop}\nolimits_{#1}}}

\theoremstyle{remark}

\def\Xint#1{\mathchoice
	{\XXint\displaystyle\textstyle{#1}}%
	{\XXint\textstyle\scriptstyle{#1}}%
	{\XXint\scriptstyle\scriptscriptstyle{#1}}%
	{\XXint\scriptscriptstyle\scriptscriptstyle{#1}}%
	\!\int}
\def\XXint#1#2#3{{\setbox0=\hbox{$#1{#2#3}{\int}$}
	\vcenter{\hbox{$#2#3$}}\kern-.5\wd0}}

\def\dashint{\Xint-}

\begin{document}

\title {Analysis on metric spaces}

\author{Mario Bonk
\affil{ Mario Bonk, Department of Mathematics, University of California, Los Angeles, {\tt mbonk@math.ucla.edu}}
\and
 Luca Capogna
\affil{ Luca Capogna, Department of Mathematical Sciences, Worcester Polytechnic Institute, {\tt lcapogna@wpi.edu}}
\and
Piotr Haj\l asz
\affil{Piotr Haj\l asz , Department of Mathematics, University of Pittsburgh, {\tt hajlasz@pitt.edu}}
\and
 Nageswari Shanmugalingam
\affil{Nageswari Shanmugalingam,  Department of Mathematical Sciences, University of Cincinnati, {\tt shanmun@uc.edu}}
\and
Jeremy Tyson
\affil{Jeremy Tyson,  Department of Mathematics, University of Illinois at Urbana-Champaign, {\tt tyson@illinois.edu}}
}

\maketitle

The subject of analysis, more specifically first-order calculus, in metric measure spaces provides a unifying framework for ideas and questions from many different fields of mathematics. One of the earliest motivations and applications of this theory arose in Mostow's work  \cite{Mostow}, in which he extended his celebrated  rigidity theorem  for hyperbolic manifolds to the  more general framework of manifolds locally modeled on negatively-curved symmetric spaces of rank one. In his proof, Mostow used the theory of quasiconformal mappings on the visual boundaries of rank-one symmetric spaces.  These visual boundaries are equipped with a sub-Riemannian structure that is locally non-Euclidean and has a fractal nature. Mostow's study of quasiconformal maps on such boundaries motivated Heinonen and Koskela \cite{HK} to axiomatize several aspects of Euclidean quasiconformal geometry in the setting of metric measure spaces, and thereby extend Mostow's work beyond the sub-Riemannian setting. The groundbreaking work \cite{HK} initiated the modern theory of analysis on metric spaces.

Analysis on metric spaces is nowadays an active and independent field, bringing together researchers from different parts of the mathematical spectrum. It has far-reaching applications to areas as diverse as geometric group theory, nonlinear PDEs, and even theoretical computer science. As a further sign of recognition, analysis on metric spaces has been included in the 2010 MSC classification as a 
category 
({\em 30L: Analysis on metric spaces}). In this short survey, we can only discuss  a small fraction of areas into which analysis on metric spaces has expanded. For more comprehensive introductions to various aspects of the subject, we invite the reader to consult the monographs \cite{Heinonenbook,HajKos, MR3363168,MR2867756,MR2662522,MR2401600,Buyalo_Schroeder_book,Hei:NonSmth}.

\vspace{0.05in}

\noindent{\bf Poincar\'e inequalities in metric spaces.}
Inspired by the fundamental theorem of calculus, Heinonen and Koskela proposed the notion of \emph{upper gradient} as a substitute for the derivative of a function on a metric measure space $(X,d,\mu)$. More precisely, $g\ge 0$ is an upper gradient for a real-valued function $u$ on $X$  if $$|u(\gamma(1))-u(\gamma(0))|\le \int_\gamma g\, ds$$ for each path $\gamma\colon [0,1]\rightarrow X$ of finite length.

 Upper gradients are not unique; but  if a function $u$ has an upper gradient $g\in L^p(\mu)$, then there is a unique {\em $p$-weak upper gradient} $g_u$ with minimal $L^p$-norm, for which the preceding inequality holds for ``almost every" curve $\gamma$.
The metric measure space $X$  is said to support a {\em $p$-Poincar\'e inequality} for some $p\ge 1$ if constants $C>0$ and $\lambda\ge 1$ exist so that for every ball $B=B(x,R)\subset X$, the inequality 
\begin{gather*}\dashint_{B}|u-u_B|d\mu 
\le CR\biggl(\,\dashint_{\lambda B}g_u^p\, d\mu\biggr)^{1/p}
\end{gather*} 
holds for all function-upper gradient pairs $(u,g_u)$. Here 
$u_B=\dashint_{B}u\, d\mu $ and ${\lambda B}=B(x,\lambda R)$.  

Over the past twenty years, many aspects of first-order  calculus 
 have been  systematically developed in the setting of 
{\em PI spaces}, that is, metric measure spaces equipped with a doubling measure and supporting a Poincar\'e inequality. For example, for PI spaces  we now have a rich theory  of Sobolev
functions which in turn lies at the foundation 
of the theory of quasconformal mappings and non-linear potential theory.  



A wealth of interesting and important examples of non-Euclidean PI spaces exist, including sub-Riemannian manifolds such as the Heisenberg group, Gromov-Hausdorff limits of manifolds with lower Ricci curvature bounds, visual boundaries of certain hyperbolic buildings, and fractal spaces that are homeomorphic to the Menger curve. 
The scope of the theory, however,   is 
not fully explored.

\vspace{0.05in}

\noindent{\bf Quasiconformal maps and nonlinear potential theory in metric spaces.}
A homeomorphism between metric spaces is said to be {\em quasiconformal} if it distorts the geometry of infinitesimal balls  in controlled fashion. Conformal maps form a special subclass for which infinitesimal balls  are mapped to infinitesimal balls.
Since the only  conformal maps between 
higher-dimensional Euclidean spaces are M\"obius transformations, quasiconformal homeomorphism  form a more flexible class 
for geometric mapping problems. For quasiconformal maps 
on PI spaces, we now have  a well-developed theory that features many of the aspects of the Euclidean theory such as Sobolev regularity, preservation of sets of measures zero, and global distortion estimates, among other things. 

A  function $u$ on a domain $\Omega$ in a metric measure space $(X,d,\mu)$ is said to be \emph{$p$-quasi-harmonic} for $p\ge 1$ if a constant $Q\ge 1$ exists so that the inequality
$$
\int_{\text{spt}\varphi}g_u^p\, d\mu\le Q \int_{\text{spt}\varphi}g_{u+\varphi}^p\, d\mu
$$
holds whenever $\varphi$ is a Lipschitz function with compact support $\text{spt}\, \varphi$ in $\Omega$. In case $Q=1$, we say that $u$ is \emph{$p$-harmonic}; this coincides with the classical Euclidean notion of a $p$-harmonic function, defined as a weak solution to the $p$-Laplace equation $$\text{div}(|\nabla u|^{p-2}\nabla u)=0.$$

Quasi-harmonic functions are useful in the study of quasiconformal mappings. For example, one can characterize quasiconformal homeomorphisms between $n$-dimensional Euclidean domains as those homeomorphisms that preserve the class of $n$-quasi-harmonic functions. A similar statement is also true for PI spaces. This generalizes the  well-known fact that planar conformal mappings are precisely the orientation-preserving homeomorphisms that preserve harmonic functions under pull-back.

The further development of potential theory in the setting of metric measure spaces leads to a classification of spaces as either $p$-parabolic or $p$-hyperbolic.  This dichotomy can be seen as a non-linear analog of the recurrence/transience dichotomy in the theory of Brownian motion. This classification is helpful in the development of a quasiconformal uniformization theory, or for a 
deeper understanding  of the links between the geometry of hyperbolic spaces and the analysis on their boundaries at infinity.

\vspace{0.05in}

\noindent{\bf Differentiability of Lipschitz functions.}
The notion of upper gradient generalizes, to metric spaces, the {\em norm} of the gradient of a $C^1$-function. It is a priori unclear how to formulate a notion of the {\em gradient} itself  (or of the differential of a function) in the absence of a linear structure.  Cheeger \cite{Cheeger} introduced a linear differential structure  for real-valued functions on metric measure spaces, and established a version of Rademacher's theorem for Lipschitz functions defined on PI spaces.
This differential structure gives rise to a finite-dimensional measurable vector bundle, the {\em generalized cotangent bundle}, over the metric space: to
each real-valued Lipschitz function $u$ corresponds an 
$L^\infty$-section $du$ of this
bundle. Moreover, the pointwise Euclidean norm $|du|$ is comparable to the minimal upper gradient
$g_u$ almost everywhere. This structure
can be used in turn to investigate second-order PDEs in divergence form,
as a basis for a theory of differential currents in metric spaces, and for many other purposes.

\vspace{0.05in}

\noindent{\bf Bi-Lipschitz embedding theorems.}
An earlier version of Rademacher's differentiation theorem, for Lipschitz maps between Carnot groups, was proved by Pansu \cite{Pansu}. Semmes observed that the Pansu--Rademacher theorem implies that nonabelian Carnot groups do not admit bi-Lipschitz copies in finite-dimensional Euclidean spaces. Moreover, such spaces do not bi-Lipschitz embed into Hilbert space, or even into any Banach space with the Radon-Nikod\'ym property (RNP). Indeed, the algebraic features of sub-Riemannian geometry have direct implications for metric questions such as bi-Lipschitz equivalence or embeddability.

The bi-Lipschitz embedding problem is intimately related to the existence of suitable differentiation theories for Lipschitz functions and maps. Roughly speaking, this relationship proceeds via incompatibility between the geometry of the cotangent bundles of the source and target spaces. In view of Cheeger's differentiation theorem, one can allow arbitrary  PI space as source spaces here 
and take  RNP Banach spaces as targets, for example. On the other hand, there is no effective differentiation theory for maps into $\ell^\infty$, because according to the Fr\'echet embedding theorem, \emph{every} separable metric space embeds isometrically into $\ell^\infty$. 

The target space $L^1$ presents an interesting intermediate case. It is not an RNP Banach space, yet deep bi-Lipschitz nonembedding theorems are available for this target. In particular, Cheeger and Kleiner \cite{ck:L1} showed that the Heisenberg group does not bi-Lipschitz embed into $L^1$. In concert with results of Lee and Naor, this fact exhibits the Heisenberg group as a geometrically natural example relevant for algorithmic questions in computer science. There has been significant additional quantitative work along these lines, culminating in Naor and Young's sharp lower bound for the integrality gap of the Goemans--Linial semidefinite program for the Sparsest Cut Problem \cite{MR3815462}. For further information, see Naor's ICM lecture \cite{naor:ICM}.

\vspace{0.05in}

\noindent{\bf Geometric measure theory on metric spaces.}
In yet another direction, Kirchheim's proof of the almost everywhere metric differentiability of Lipschitz mappings into metric spaces led to far-reaching generalizations of the area and co-area formulas.
Subsequently, Ambrosio and Kirchheim \cite{Ambrosio-Kirchheim} developed an extension of the Federer-Fleming theory of currents in complete metric spaces, thus opening
a new chapter in geometric measure theory leading to a 
study of 
(quantitative) rectifiability in metric spaces. This notion 
has been
further
developed  
in sub-Riemannian settings, especially in the Heisenberg group, but many questions remain open. 
These theories have been relevant, for example,  in  the work of Naor and Young \cite{MR3815462} mentioned above, where quantitative rectifiability of surfaces in the Heisenberg group  is prominently featured.

\vspace{0.05in}

\noindent{\bf Dynamics and analysis on metric spaces.}
An interesting source of examples of spaces that can be studied with the methods of quasiconformal geometry and analysis on metric spaces are fractals that arise from self-similar  or dynamical constructions such as limit sets of Kleinian  groups, Julia sets of rational maps, or attractors  of iterated function systems. Often the geometry of these spaces is too ``rough" to expect finer analytic properties such the Poincar\'e inequality to hold. However, if these spaces admit a good first-order calculus, then striking consequences often emerge. For example, Cannon's well-known conjecture in geometric group theory predicts that a Gromov hyperbolic group $G$ admits a geometric action on hyperbolic $3$-space if its boundary at infinity $\partial_\infty G$ is a topological $2$-sphere.  While the conjecture is still open, one can show that the desired conclusion is true if $\partial_\infty G$ (equipped with a  visual metric) has good analytic properties, say, if it is 
quasisymmetrically equivalent to 
a PI space. For more information see the  ICM lectures \cite{bonk:ICM} and \cite{kleiner:ICM}.

The problem of deciding when a metric space is quasisymmetrically equivalent to a space with ``better" analytic properties can be seen as a generalization
of classical uniformization theorems in complex analysis, at least for low-dimensional fractals such as  Sierpi\'nski carpets or fractal $2$-spheres. In order to study such problems, one often employs concepts from classical complex analysis in a metric space setting. For example, the modulus of a path family, originally introduced by Ahlfors and Beurling in the complex plane, now plays a prominent role in much recent work on mapping theory in general, abstract metric spaces.

\vspace{0.05in}

\noindent{\bf Conclusion.}  This brief note barely hints at the breadth and the depth of the problems of current concern in the theory of analysis on metric spaces. The 2019 AMS MRC {\it Analysis on Metric Spaces} will address a number of questions that have been the subject of much recent investigation, but are far from being completely understood.

\vspace{0.05in}

\noindent{\bf Acknowledgements.} The second author is supported by the Simons Foundation, and the other authors by the National Science Foundation.

\bibliography{general-BIBLIO-2}
\bibliographystyle{amsalpha}

\end{document}